\documentclass[12pt]{article}
\usepackage{amsmath,amsfonts,amsthm,amscd,amssymb,graphicx}

\usepackage{graphicx}

\usepackage{amssymb}
\usepackage[colorlinks,citecolor=blue,linkcolor=blue]{hyperref}
\usepackage{xcolor}

\def\eps{\varepsilon}

\newtheorem{theorem}{Theorem}[section]
\newtheorem{lemma}[theorem]{Lemma}
\newtheorem{e-proposition}[theorem]{Proposition}

\newtheorem{e-definition}[theorem]{Definition\rm}

\def\og{\leavevmode\raise.3ex\hbox{$\scriptscriptstyle\langle\!\langle$~}}
\def\fg{\leavevmode\raise.3ex\hbox{~$\!\scriptscriptstyle\,\rangle\!\rangle$}}

\def\beq{\begin{equation}}
\def\eeq{\end{equation}}

\begin{document}



%

\centerline{\Large \bf Long waves instabilities}

\bigskip

\centerline{D. Bian\footnote{Beijing Institute of Technology, School of Mathematics and Statistics, Beijing, China}, 
E. Grenier\footnote{UMPA, CNRS  UMR $5669$,Ecole Normale Sup\'erieure de Lyon, Lyon, France}}



\subsubsection*{Abstract}


The aim of this paper is to give a detailed presentation of long wave instabilities of shear layers for
Navier Stokes equations, and in particular to give a simple and easy to read presentation of the study of Orr Sommerfeld equation
and to detail the analysis of its adjoint. Using these analyses we prove the existence of long wave instabilities in the case
of slowly rotating fluids, slightly compressible fluids and for Navier boundary conditions, under  smallness conditions.


\section{Introduction}


Let us first consider incompressible Navier Stokes equations in an half space
 \beq \label{NS1} 
\partial_t u^\nu + (u^\nu   \cdot \nabla) u^\nu - \nu \Delta u^\nu + \nabla p^\nu = f^\nu,
\eeq
\beq \label{NS2}
\nabla \cdot u^\nu = 0,
\eeq
together with the Dirichlet boundary condition
\beq \label{NS3} 
u^\nu = 0 \qquad \hbox{for} \qquad z = 0.
\eeq
The study of the inviscid limit $\nu \to 0$ of this system has been extensively studied (see for instance
\cite{GR08,DGV2,GGN3, Mae,SammartinoCaflisch1,SammartinoCaflisch2} and the references therein). 
In this paper we are interested in the stability of a shear layer profile
$U(z) = (U_s(z), 0)$. Note that this shear layer profile is a stationary solution of Navier Stokes equations
provided we add the forcing term $f^\nu = (- \nu \Delta U_s, 0)$. 
The stability of such shear layer profiles has been intensively studied in physics since the beginning of the
twentieth century, see \cite{Blasius,Reid,Hei,Lin0,LinBook,Schlichting} and the references therein.  
We recall that such shear layers are unstable in $3$ dimensions if and only if they are unstable in $2$ dimensions
(so called Squire's theorem).
We also recall that non trivial shear layer profiles are always linearly unstable with respect to Navier Stokes equations.
Two classes of instabilities appear as follows:

\begin{itemize}

\item "Inviscid instabilities": instabilities which persist as $\nu$ goes to $0$. According to Rayleigh's criterium,
such instabilities only occur for profiles
$U_s$ with inflection points, and exhibit scales in $t$, $x$ and $z$ of order $O(1)$.
This first kind of instabilities is now well understood, both at the linear and nonlinear levels \cite{GN, GN1}.

\item "Long wave instabilities": these instabilities arise even in the case of concave profiles $U_s$, such that
$U_s'' < 0$. They
are characterized by a strong spatial anisotropy since their sizes are of order $O(1)$ in $z$ but of order
$O(\nu^{-1/4})$ in $x$. Moreover they grow very slowly, within time scales of order $O(\nu^{-1/2})$.
This second type of instabilities is more delicate to study, even at the linear level.

\end{itemize}

The purpose of this paper is to give an educational presentation of this second kind of instabilities, by recalling
recent results on Orr Sommerfeld equations, and to give the first mathematical study of the adjoint of Orr Sommerfeld
equation. We then apply these studies to systems close to incompressible Navier Stokes equations with Dirichlet 
boundary condition. We in particular formally prove that there exists such long wave instabilities for
incompressible Navier Stokes equations with Navier boundary conditions, for slightly rotating fluids and
for compressible Navier Stokes equations in the low Mach number regime, all under a smallness assumption.
We only give formal proofs. Rigorous approaches are straightforward using the methods used in \cite{GN}. 
This present work is a first step in a general program to study the nonlinear instability in the viscous boundary layer
\cite{Bian2}.


\section{Orr Sommerfeld equations}


We first recall the main results on Orr Sommerfeld equations. The corresponding proofs may be found for instance
in \cite{GN1}.


\subsection{Introduction}


Let $L$ be the linearized Navier Stokes operator near the shear layer profile $U$, namely
\beq \label{linearNS}
L u = (U  \cdot \nabla) u + (u \cdot \nabla) U - \nu \Delta u + \nabla q,
\eeq
with $\nabla  \cdot u = 0$ and Dirichlet boundary condition. We want to study the resolvent of $L$, namely to
study the equation
\beq \label{resolvant}
(L - \lambda) u = f
\eeq
where $f$ is a given forcing term.
Following the classical analysis we take advantage of the incompressibility condition to introduce the stream function,
take its Fourier transform in $x$ and its Laplace transform in time, and thus look for velocities of the form
$$
u = \nabla^\perp \Bigl( e^{i \alpha (x - c t) } \psi(z) \Bigr)=e^{i \alpha (x - c t) }(\partial_z\psi, -i\alpha \psi).
$$
Note that  $\lambda = i \alpha c$. We also take the Fourier and Laplace transform of the forcing term $f$
$$
f = \Bigl( f_1(z),f_2(z) \Bigr) e^{i \alpha (x - c t) } .
$$
Taking the curl of (\ref{resolvant}) we then get
\beq \label{Orrmod}
Orr_{c,\alpha,\nu}(\psi) :=  (U_s - c)  (\partial_z^2 - \alpha^2) \psi 
 - U_s''  \psi  
- { \nu \over i \alpha}   (\partial_z^2 - \alpha^2)^2 \psi =- {\nabla \times f \over i \alpha},
\eeq
where
$$
\nabla \times (f_1,f_2) = i \alpha f_2 - \partial_z f_1.
$$
The Dirichlet boundary condition
gives 
\beq \label{condOrr}
\psi(0) = \partial_z \psi(0) = 0.
\eeq
Let
\beq \label{epsilon}
\eps = {\nu \over i \alpha} .
\eeq
As $\nu$ goes to $0$, the $Orr$ operator degenerates into the classical Rayleigh operator
\beq \label{Rayleigh}
Ray_\alpha(\psi) = (U_s - c)  (\partial_z^2 - \alpha^2) \psi 
 - U_s''  \psi
 \eeq
 which is a second order operator, together with the boundary condition 
 \beq \label{condRay}
 \psi(0) = 0.
 \eeq
 Two cases arise, depending on whether this Rayleigh operator has unstable eigenvalues or not.
 In the first case there exists $\alpha$, $c$ and $\psi_R$, which satisfy \eqref{condRay}, 
 such that 
 $$
 Ray_\alpha(\psi_R) = 0.
 $$
 Starting from this solution, through a perturbative analysis, it is possible to construct eigenmodes $\psi_\nu$
 of the Orr Sommerfeld operator  with corresponding eigenvalues $c_\nu$ which satisfy (\ref{condOrr}) and
$Orr(\psi_\nu) = 0$, and it is also possible to prove the nonlinear instability of the corresponding shear layer profiles. 
This case is now well understood and  we in particular refer to \cite{GN} and \cite{GN1} for further details.
 
We will thus focus on the second case, where Rayleigh's operator has no unstable eigenvalue.

 
 \subsection{Rayleigh equation}
 

 Let us first focus on Rayleigh's equation. It is a second order differential equation. For $c$ close to
 zero, we note that $U_s - c$ almost vanishes.
 More precisely for small $c$, let $z_c$ be such that 
 \beq \label{yc}
 U_s(z_c) = c .
 \eeq
 Such  a complex number $z_c$ is called a "critical layer". At $z_c$, the Rayleigh operator degenerates from 
 a second order operator to an operator of order $0$ and thus has a "turning point" at $z = z_c$.

We therefore expect that there exist two independent solutions of Rayleigh operator, one regular at $z_c$ and the
other one which is singular near this point. It turns out that it is possible to get a very accurate description
of these two independent solutions for small $\alpha$.

Namely we first observe that for $\alpha = 0$, Rayleigh operator degenerates into
 $$
 Ray_0(\psi) = (U_s - c) \partial_z^2 \psi - U_s'' \psi.
 $$
 In particular
 \beq \label{psimoins0}
  \psi_{-,0}(z) = U_s(z) - c
 \eeq
 is an explicit solution of this limiting operator. 
 An independent solution $\psi_{+,0}$ can be computed using the Wronskian
 $$
 W(\psi_{-,0},\psi_{+,0}) = 1.
 $$ 
 More precisely,
 \beq \label{psiplus0}
 \psi_{+,0}(z) = C(z) \psi_{-,0}(z)
 \eeq
 where
 $$
 C'(z) = (U_s(z) - c)^{-2} .
 $$
 This other solution behaves linearly at infinity and has a $(z - z_c) \log(z - z_c)$ singularity at $z = z_c$.
  
 Using a perturbative argument it is possible to prove that, for $\alpha$ small enough, 
 the Rayleigh operator has two independent solutions, one, $\psi_{-,\alpha}$ which goes to
 $0$ at infinity and is close to $U_s - c$ for bounded $z$, and another
 one, $\psi_{+,\alpha}$, which diverges at $+ \infty$ and which has a $(z - z_c) \log(z - z_c)$ singularity at $z_c$.

\subsection{Solutions of Orr Sommerfeld \label{sec23}}

Let us go back to the Orr Sommerfeld operator. 
Note that $Orr$ is a fourth order differential operator. It can be proven that there exists four independent solutions
to this operator (not taking into account any boundary condition),
two with a "slow" behavior, called $\psi_{s,\pm}$, and two with a "fast" behavior, called $\psi_{f,\pm}$.
The "-" subscript refers to solutions which go to $0$ at infinity and the "+" subscript to solutions diverging as $y$ goes
to $+ \infty$.

The "slow" solutions have bounded derivatives. For such solutions the $Orr$ operator reduces to Rayleigh's operator
as $\nu$ goes to $0$. More precisely,
$$
Orr = Ray - \eps Diff,
$$
where
$$
Diff(\psi) =   (\partial_z^2 - \alpha^2)^2 \psi .
$$
For "slow" solutions, $Diff$ may be treated as a small perturbation. For $\psi_{-,\alpha}$ we directly have
$Orr(\psi_{-,\alpha}) = O(\eps)$. For $\psi_{+,\alpha}$ the situation is more delicate since
$Orr(\psi_{+,\alpha})$ may be large near $z_c$. However, using an iterative scheme, treating $Diff$ as a perturbation, it is possible to construct genuine solutions
of the $Orr$ operator, which are close to $\psi_{\pm,\alpha}$.
In particular it can be proven that $\psi_{s,-}$ satisfies
\beq \label{psism}
\psi_{s,-}(0) = - c + \alpha {U_+^2 \over U_s'(0)} + O(\alpha^2),
\eeq
\beq \label{psisd}
\partial_z \psi_{s,-}(0) = U'_s(0) + O(\alpha),
\eeq
\beq \label{psisdd}
\partial_z^2 \psi_{s,-}(0) = O(1),
\eeq
where $U_+=\lim_{z\rightarrow +\infty}U_s(z)$.

On the contrary "fast" solutions have very large gradients and higher derivatives. For these solutions 
Orr Sommerfeld may be approximated by its higher order derivatives, namely by
\beq \label{modAiry}
(U_s - c) \partial_z^2 - \eps \partial_z^4 = Airy  \, \partial_z^2,
\eeq
where $Airy$ is the modified Airy operator 
$$
Airy =  (U_s - c) - \eps \partial_z^2.
$$
More precisely,
$$
Orr = Airy  \,  \partial_z^2 + Rem,
$$
where the remainder operator $Rem$ contains derivatives of lower order, and is therefore smaller with respect
to $Airy \, \partial_z^2$ for "fast" solutions.

The next step is to construct solutions of $Airy \, \, \partial_z^2 = 0$, or, up to two integrations, to construction
solutions of $Airy = 0$. Near $z_c$, this operator may be approximated by the classical Airy operator
\beq \label{linAiry}
U_s'(z_c) (z - z_c) - \eps \partial_z^2 .
\eeq
Solutions of (\ref{linAiry}) may be explicitly expressed as combinations of the classical Airy's functions
$Ai$ and $Bi$, which are solutions of the classical Airy equation $z \psi = \partial_y^2 \psi$.
More precisely, $Ai(\gamma (z - z_c))$ and $Bi(\gamma (z - z_c))$ are solutions of (\ref{linAiry}) provided
we choose
\beq \label{gamma}
\gamma = \Bigl( {i \alpha U_s'(z_c) \over \nu} \Bigr)^{1/3}.
\eeq
Now solutions of $Airy = 0$ may be constructed, starting from solutions of (\ref{linAiry}) through the
so called Langer's transformation (which is a transformation of the phase and amplitude of the solution).
To go back to (\ref{modAiry}) we then have to integrate twice these solutions.
As a consequence, fast solutions to the Orr Sommerfeld equation may be expressed in terms of
second primitives of the classical Airy functions $Ai$ and $Bi$.
 Let us call $Ai(2,.)$ and $Bi(2,.)$ this second primitives of $Ai$ and $Bi$.
 
It can be proven that $\psi_{f,-}$ satisfies
 \beq \label{psif}
\psi_{f,-}(0) = Ai(2,-\gamma z_c) + O(\alpha),
\eeq
\beq \label{psifd}
\partial_z \psi_{f,-}(0) = \gamma Ai(1, - \gamma z_c) + O(1),
\eeq
\beq \label{psifdd}
\partial_z^2 \psi_{f,-}(0) = \gamma^2 Ai(- \gamma z_c) + O(\alpha^{-1}).
\eeq

\subsection{Dispersion relation}

An eigenmode of Orr Sommerfeld equation is a combination of these particular solutions which goes to $0$ at infinity
and which vanishes, together with its first derivative, at $z  = 0$. As an eigenmode must go to $0$ as $y$ goes to 
infinity, it is a linear combination of $\psi_{f,-}$ and $\psi_{s,-}$ only.
There should exist nonzero constants $a$ and $b$ such that 
$$
a \psi_{f,-} (0) + b \psi_{s,-} (0) = 0
$$ and
$$
a \partial_z  \psi_{f,-} (0) + b \partial_z \psi_{s,-} (0) = 0.
$$
Equivalently the determinant
$$
E = \left| \begin{array}{cc}
 \psi_{f,-} (0) &  \psi_{s,-} (0) \cr
  \partial_z  \psi_{f,-} (0) &  \partial_z  \psi_{s,-} (0) \cr
\end{array} \right|
$$
must vanish.
The dispersion relation is therefore
\beq \label{disper}
{   \psi_{f,-} (0) \over \partial_z  \psi_{f,-} (0)}
= {  \psi_{s,-} (0) \over \partial_z  \psi_{s,-} (0)}
\eeq
or
\beq \label{disper2}
\alpha {U_+^2 \over U_s'(0)^2}  - {c \over U_s'(0)} = \gamma^{-1} {Ai(2, - \gamma z_c) \over Ai(1,-\gamma  z_c)} 
+O(\alpha^2).
\eeq
We will focus on the particular case where $\alpha$ and $c$ are both of order $\nu^{1/4}$. It turns out 
that this is an area where instabilities occur, and we conjecture that it is in this region that the most
unstable instabilities may be found.
We rescale $\alpha$ and $c$ by $\nu^{1/4}$ and introduce
$$
\alpha = \alpha_0 \nu^{1/4}, \quad 
c = c_0 \nu^{1/4}, \quad
Z = \gamma z_c,
$$
which leads to
\beq \label{disper3}
\alpha_0 {U_+^2 \over U_s'(0)^2}  - {c_0 \over U_s'(0)} = 
 {1 \over (i \alpha_0 U'_s(z_c))^{1/3}}  {Ai(2, - Z) \over Ai(1,- Z)} + O(\nu^{1/4}).
\eeq
Note that as $U_s(z_c) = c$, 
$$
z_c = U_s'(0)^{-1} c + O(c)
$$
and
$$
Z =  \Bigl(i  U_s'(z_c) \Bigr)^{1/3} \alpha_0^{1/3} U_s'(0)^{-1} c_0 + O(\nu^{1/4}) .
$$
Note that the argument of $Z$ equals $\pi / 6$.
We then introduce the following function, called Tietjens function, of the real variable $z$
$$
Ti(z) = {Ai(2,z e^{- 5 i \pi / 6}) \over z e^{- 5 i \pi / 6} Ai(1,z e^{- 5 i \pi / 6})} .
$$
At first order the dispersion relation becames
\beq \label{dispersionlimit}
\alpha_0 {U_+^2 \over U_s'(0)} = c_0 \Bigl[ 1 - Ti(- Z e^{5 i \pi / 6} ) \Bigr] .
\eeq
This dispersion relation can be numerically investigated. Note that it only depends on the limit $U_+$
of the horizontal velocity at infinity and on $U_s'(0)$. It can be proven that provided $\alpha_0$ is large enough
there exists $c_0$ with $\Im c_0 > 0$, leading to an unstable mode.

\begin{figure}
\includegraphics[width = 10cm]{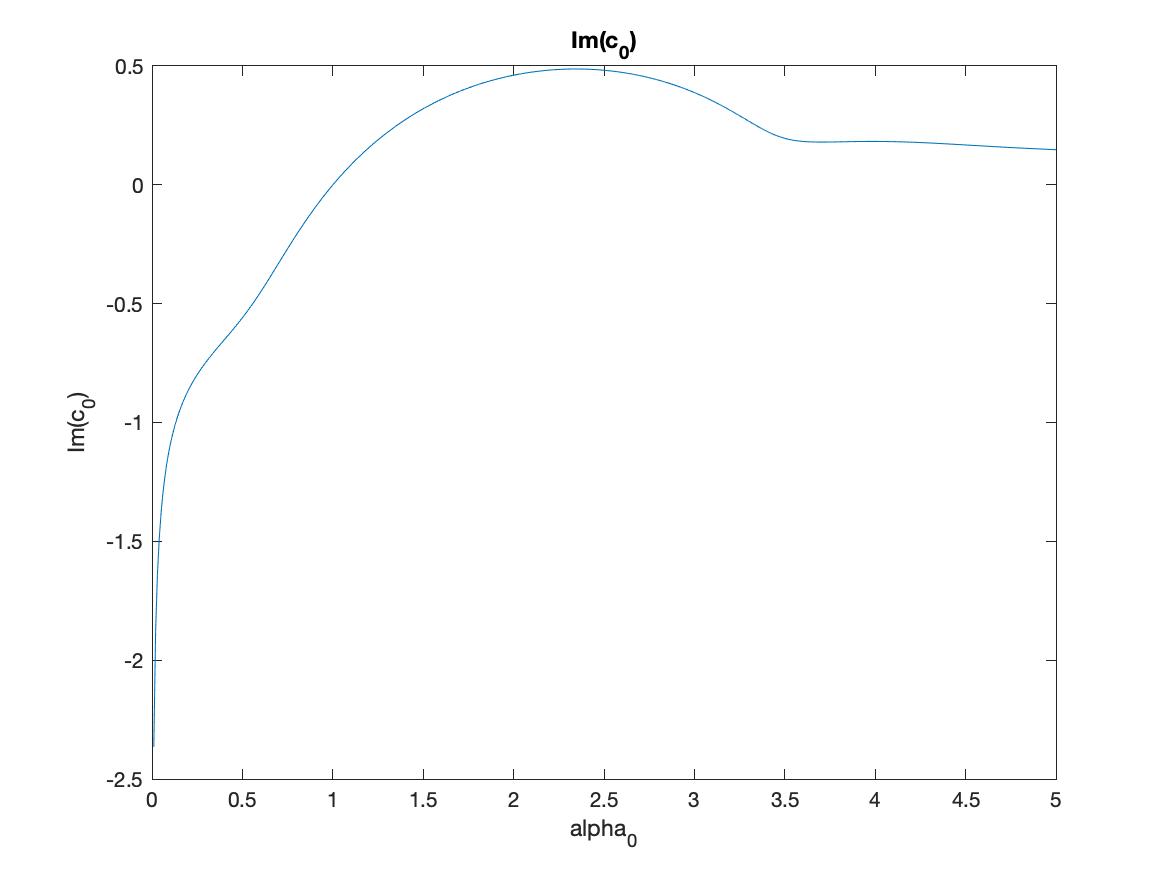}
\caption{$\Im c_0$ as a function of $\alpha_0$}
\label{valeurc}
\end{figure}

\begin{figure}
\includegraphics[width = 10cm]{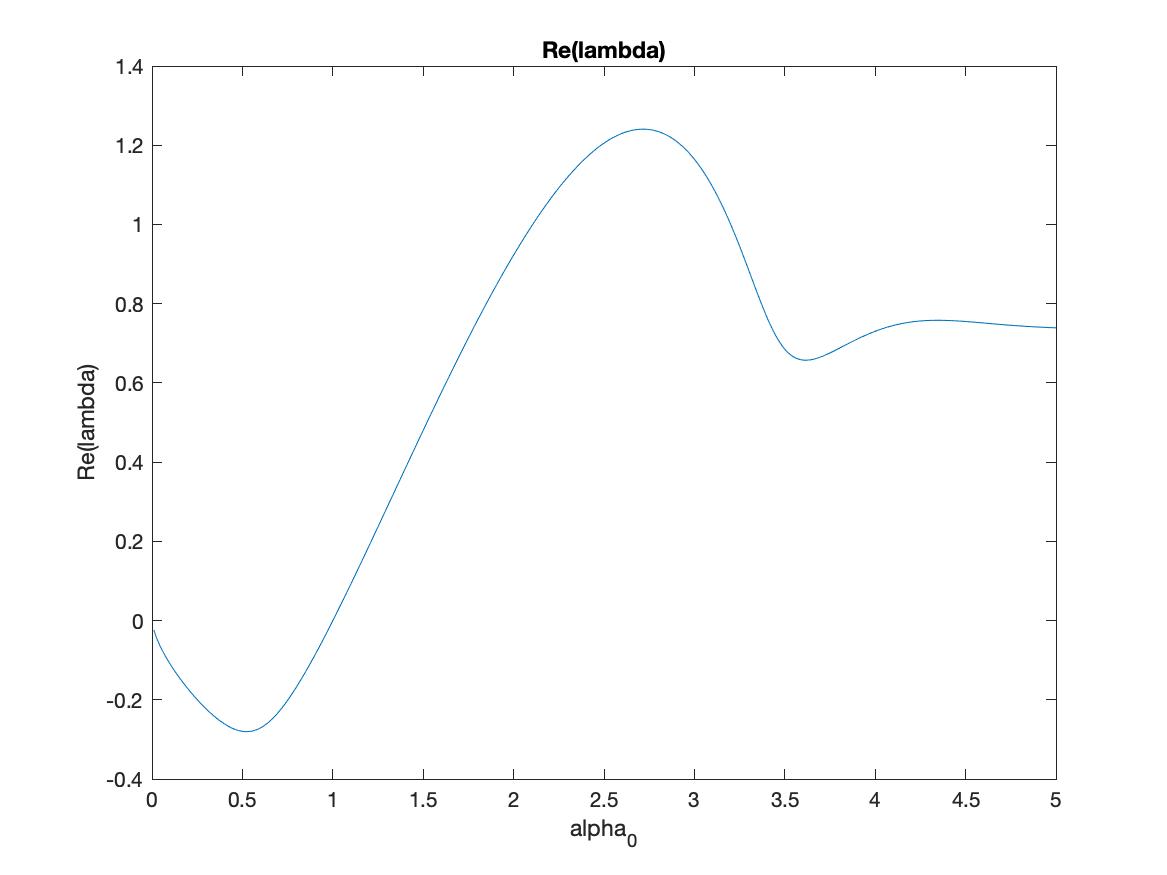}
\caption{$\Re \lambda$ as a function of $\alpha_0$}
\label{valeurlambda}
\end{figure}

Let us numerically illustrate this dispersion relation in the case $U_s'(0) = U_+ = 1$.
Figure \ref{valeurc} shows the imaginary part of $c_0$ as a function of $\alpha_0$. We see that there exists 
a constant $\alpha_c$ such that $\Im c_0 > 0$ if $\alpha_0 > \alpha_c$ and $\Im c_0 < 0$ if
$\alpha_0 < \alpha_c$.
Figure \ref{valeurlambda} shows $\Re \lambda = \alpha_0 \Im c_0$ as a function of $\alpha_0$.
We see that there exists a unique global maximum to this function, at $\alpha_0 = \alpha_M \sim 2.8$.

\subsection{Description of the linear instability}


Let us now detail the linear instability. Its stream function $\psi_{lin}$ is of the form
$$
\psi_{lin} = b \psi_{s,-} + a \psi_{f,-} .
$$
Choosing $b = 1$ we see that $a = O(\nu^{1/4})$, hence
\beq \label{psilin}
\psi_{lin}(z) = U_s(z) - c + \alpha {U_+^2 \over U_s'(0)} + a Ai \Bigl(2, \gamma (z - z_c) \Bigr) +  O(\nu^{1/2} ).
\eeq
The corresponding horizontal and vertical velocities $u_{lin}$ and $v_{lin}$ are given by
\beq \label{vlin}
u_{lin} = \partial_z \psi_{lin} = U_s'(z) +\gamma a Ai \Bigl(1, \gamma (z - z_c) \Bigr) + O(\nu^{1/4}),
\eeq
and 
\beq \label{ulin} 
v_{lin} = -i\alpha \psi_{lin} = O(\nu^{1/4}).
\eeq
Note that $\gamma a = O(1)$. The first term in (\ref{vlin}) may be seen as a "displacement velocity", corresponding
to a translation of $U_s$. The second term is of order $O(1)$ and located in the boundary layer, namely
within a distance $O(\nu^{1/4})$ to the boundary. Note that the vorticity
	\beq \label{check-use-number}
\omega_{lin} = -(\partial_z^2 - \alpha^2) \psi_{lin}(z) 
= -U_s''(z) - \gamma^2 a Ai \Bigl( \gamma (z - z_c) \Bigr) 
\eeq

is large in the critical layer (of order $O(\nu^{-1/4})$).


\section{Adjoint Orr Sommerfeld operator}


Let us now turn to the study of the adjoint of Orr Sommerfeld operator.


\subsection{Definition of the adjoint}


We will consider the classical $L^2$ product between two stream functions $\psi_1$ and $\psi_2$, namely
$$
(\psi_1,\psi_2) = \int \psi_1 \bar \psi_2 dx .
$$
We have
$$
(Orr (\psi_1), \psi_2) = (\psi_1, TOrr (\psi_2)),
$$
where
$$
TOrr (\psi) =  (\partial_z^2 - \alpha^2) (U_s - \bar c)   \psi - U_s''  \psi  
+ { \nu \over i \alpha}   (\partial_z^2 - \alpha^2)^2 \psi,
$$
Taking the complex conjugate we define the adjoint of Orr Sommerfeld equation to be
\beq \label{Orrad}
Orr^t_{c,\alpha,\nu}(\psi) :=  (\partial_z^2 - \alpha^2) (U_s - c)   \psi 
- U_s''  \psi  
- { \nu \over i \alpha}   (\partial_z^2 - \alpha^2)^2 \psi,
\eeq
with boundary conditions $\psi(0) = \partial_z \psi(0) = 0$.
We also introduce the corresponding adjoint of the Rayleigh operator
\beq \label{Rayad}
Ray^t =  (\partial_z^2 - \alpha^2) (U_s - c)   \psi 
- U_s''  \psi  
\eeq
and the adjoint Airy operator
\beq \label{Airyad}
{\cal A} = (\partial_z^2 - \alpha^2) \, Airy .
\eeq
We already know that the spectrum of Orr Sommerfeld operator and that of its adjoint is the same. 
Let $\psi_1$ be an eigenvector of $Orr$ with corresponding eigenvalue $\lambda_1$ and
let $\psi_2$ be an eigenvector of $Orr^t$ with corresponding eigenvalue $\lambda_2$.
Then, multiplying $Orr(\psi_1) = \lambda_1 \psi_1$ by $\bar \psi_2$, $Orr^t(\psi_2) = \lambda_2 \psi_2$ by $\bar \psi_1$ and
combining we get (see \cite{Schensted}, chapter $2$)
\beq \label{ortho}
(\lambda_1 - \lambda_2) \int \psi_1 (\partial_z^2 - \alpha^2) \bar \psi_2 dx = 0.
\eeq
We can therefore normalize the eigenvectors $\psi_1$ and $\psi_2$ such that
\beq \label{ortho2} 
\int \psi_1 (\partial_z^2 - \alpha^2) \bar \psi_2 dx = \delta_{\lambda_1 = \lambda_2}.
\eeq
Let $v_1 = \nabla^\perp \psi_1$ and $v_2 = \nabla^\perp \psi_2$, then we get
\beq \label{ortho3}
\int v_1 \bar v_2 dx =  \delta_{\lambda_1 = \lambda_2}.
\eeq


\subsection{Study of the adjoint Rayleigh operator}


The main observation is that $Ray$ and $Ray^t$ are conjugated through $U_s - c$, namely
\beq \label{conjugation}
Ray^t (\psi) =  (U_s - c)^{-1} Ray \Bigl( ( U_s - c) \psi \Bigr) .
\eeq
As a consequence
\beq \label{conjugation2}
(Ray^t)^{-1}(f ) =   (U_s - c)^{-1} Ray^{-1} \Bigl( ( U_s - c) f \Bigr) .
\eeq
All the study of $Ray^t$ can be thus deduced from that of $Ray$, as developed in \cite{GN}, a reference that we now closely follow.


\subsubsection{The case $\alpha = 0$}


For $\alpha = 0$, the adjoint of Rayleigh operator reduces to
$$
Ray^t_0(\psi) = \partial_z^2 \Bigl[ (U_s - c) \psi \Bigr] - U''_s \psi .
$$
Using (\ref{conjugation}), there exists two particular solutions to $Ray_0^t(\psi) = 0$, namely 
\beq \label{solRayt}
\psi_{-,0}^t = 1, \qquad \psi_{+,0}^t = {\psi_{+,0} \over U_s(z) - c} ,
\eeq
where $\psi_{+,0}$ is a growing solution to $Ray_0 = 0$.
Note in particular that $\psi_{+,0}^t$ has a singularity of the form $(U_s(z) - c)^{-1}$, and is 
therefore much more singular than $\psi_{+,0}$, which is bounded near the critical layer.
Note that
$$
W \Bigl[ \psi^t_{-,0}, \psi^t_{+,0} \Bigr] = {1 \over (U_s - c)^2}.
$$
The Green function of $Ray_0^t$ can be deduced from that of $Ray_0$ and is (see \cite{GN})
$$
G_{R,0}(x,z) = \left\{ \begin{array}{rrr} - (U_s(z)-c)^{-1} \psi_{-,0}(z) \psi_{+,0}(x), 
\quad \mbox{if}\quad z>x,\\
- (U_s(z)-c)^{-1} \psi_{-,0}(x) \psi_{+,0}(z), \quad \mbox{if}\quad z<x.\end{array}\right.
$$
%
%
Note that this Green function is singular in $z$, with a $(U(z) - c)^{-1}$ singularity.
The inverse of $Ray_0^t$ is explicitly given by
\begin{equation}\label{def-RayS0}
RaySolver_0^t(f) (z) : =  \int_0^{+\infty} G_{R,0}(x,z) f(x) dx.
\end{equation}
As in \cite{GN}, we define
$$
X^\eta = L_\eta^{\infty} = \Bigl\{ f \quad | \quad \sup_{z \ge 0} | f(z) | e^{\eta z} < + \infty \Bigr\}.
$$
and the spaces $Y^\eta$ defined by: a function $f$ lies in $Y^\eta$ if for any $z \ge 1$,
$$
|f(z)| + |\partial_z f(z) | + | \partial_z^2 f (z) |  \le C e^{-\eta z}
$$
and if for $z \le 1$,
$$
| \partial_z f(z) | \le C (1 + | \log (z - z_c) |  ) , 
$$
and
$$| \partial_z^2 f(z) | \le C (1 + | z - z_c |^{-1} ).
$$
The best constant $C$ in the previous bounds defines the norm $\| f \|_{Y^{\eta}}$. 
Note that, following Lemma 3.3 of \cite{GN}, if $(U_s(z) - c) f \in X^\eta$ 
then there exists a solution $\phi$ to $Ray_0^t(\phi) = f$
such that $(U_s(z) - c) \phi \in Y^0$, and
$$
\| (U_s(z) - c) \phi \|_{Y^0} \le C ( 1 + | \log \Im c | ) \| (U_s(z) - c) f \|_{X^\eta}.
$$
Again, note the singularity at $z = z_c$.


\subsubsection{Approximate Green function when $\alpha \ll1$}


We again follow \cite{GN}.
Let $\psi^t_{-,0}$ and $\psi^t_{+,0}$ be the two solutions of $Ray^t_0(\psi) = 0$ that are constructed above. 
We now construct an approximate Green function to the adjoint Rayleigh equation for $\alpha > 0$.
To proceed, let us introduce
\begin{equation}\label{def-phia12}
	\phi_{-,\alpha } = \psi^t_{-,0} e^{-\alpha z} = e^{-\alpha z}  ,\qquad \phi_{+,\alpha} = \psi^t_{+,0} e^{-\alpha z}.
\end{equation}
A direct computation shows that their Wronskian determinant equals to
$$
W[\phi_{-,\alpha},\phi_{+,\alpha}] =  \partial_z \phi_{+,\alpha} \phi_{-,\alpha} - \phi_{+,\alpha} \partial_z \phi_{-,\alpha} 
$$
$$
= e^{-2\alpha z}W[\psi^t_{-,0},\psi^t_{+,0}] = {e^{-2 \alpha z} \over (U_s(z) - c)^2}.
$$ 
Note that the Wronskian vanishes at infinity since both functions have the same behavior at infinity, and is singular at $z = z_c$. 
In addition, 
\begin{equation}\label{Ray-phia12}
Ray^t_\alpha(\phi_{\pm,\alpha}) = - 2 \alpha \partial_z \psi_{\pm,0} e^{-\alpha z}.
\end{equation}
We are now led to introduce an approximate Green function $G_{R,\alpha}(x,z)$, which is defined as follows
$$
G_{R,\alpha}^t(x,z) = \left\{ \begin{array}{rrr} - (U_s(z)-c)^{-1} e^{-\alpha (z-x)}  \psi_{-,0}(z) \psi_{+,0}(x),
 \quad \mbox{if}\quad z>x\\
- (U_s(z)-c)^{-1} e^{-\alpha (z-x)}  \psi_{-,0}(x) \psi_{+,0}(z), \quad \mbox{if}\quad z< x.\end{array}\right.
$$
Note that
$$
G^t_{R,\alpha} = { U_s(x) - c \over U_s(z) - c} G_{R,\alpha} .
$$
Similar to $G_{R,0}^t(x,z)$, the Green function $G_{R,\alpha}^t(x,z)$ is ``singular'' near $z_c$.
By a view of (\ref{Ray-phia12}), 
\begin{equation}\label{id-Gxz}
	Ray^t_\alpha (G^t_{R,\alpha}(x,z)) = \delta_{x}  + E_{R,\alpha}^t(x,z),
\end{equation}
for each fixed $x$, where the error kernel $E_{R,\alpha}^t(x,z)$ is defined by  
$$
E_{R,\alpha}^t(x,z) = 
\left\{ \begin{array}{rrr}
- 2 \alpha  (U_s(z)-c)^{-1} e^{-\alpha (z-x)} \partial_z \psi_{-,0}(z) \psi_{+,0}(x), \quad \mbox{if}\quad z>x\\
- 2 \alpha (U_s(z)-c)^{-1}e^{-\alpha (z-x)}\  \psi_{-,0}(x) \partial_z \psi_{+,0}(z), \quad \mbox{if}\quad z< x.\end{array}\right.
$$
Then an approximate inverse of the operator $Ray_\alpha$ can be defined by
\begin{equation}\label{def-RaySa}
	RaySolver_\alpha^t(f)(z) 
	:= \int_0^{+\infty} G_{R,\alpha}^t(x,z) f(x) dx,
\end{equation}
and the related error operator takes the form of
\begin{equation}\label{def-ErrR}
	Err_{R,\alpha}^t(f)(z) :=  \int_0^{+\infty} E_{R,\alpha}^t(x,z) f(x) dx.
\end{equation}
Note that
$$
RaySolver_\alpha^t (f) = (U_s - c)^{-1} RaySolver_\alpha \Bigl( (U_s - c) f \Bigr) ,
$$
and similarly for $E_{R,\alpha}^t$.
Using (\ref{conjugation2}) and Lemma $3.4$ of \cite{GN}, we get the following Lemma.
\begin{lemma}\label{lem-RaySa} 
	Assume that $\Im c > 0$. For any $f\in {X^{\eta}}$,  with $\alpha<\eta$, 
	the function $RaySolver_\alpha^t(f)$ is well-defined in $Y^{\alpha}$, and satisfies 
	$$ 
	Ray_\alpha^t(RaySolver^t_\alpha(f)) = f + Err^t_{R,\alpha}(f).
	$$
	Furthermore, there hold  
	\begin{equation}\label{est-RaySa}
	\| (U_s - c) RaySolver_\alpha^t(f)\|_{Y^{\alpha}} \le C (1+|\log \Im c|) \| (U_s - c) f\|_{{X^{\eta}}},
	\end{equation}
	and 
	\begin{equation}\label{est-ErrRa} 
	\| (U_s - c) Err_{R,\alpha}^t(f)\|_{Y^{\eta}} \le C | \alpha |   (1+|\log (\Im c)|)  \| (U_s - c) f\|_{X^{\eta}} ,
	\end{equation} 
	for some universal constant $C$. 
\end{lemma}

%


\subsection{Transposed Airy operator and construction of $\phi^{t,app}_{f,-}$}


Near the critical layer we will have to study 
$$
{\cal A} = (\partial_z^2 - \alpha^2) Airy
$$
where
$$
Airy = (U_s - c ) \psi - \eps \partial_z^2 \psi .
$$
We have
$$
Orr^t = {\cal A} - U_s'' \psi.
$$
Note that ${\cal A}$ is easily inverted into
$$
{\cal A}^{-1} = Airy^{-1} (\partial_z^2 - \alpha^2)^{-1}  .
$$
The study of $Airy$ and $Airy^{-1}$ has been fully detailled in \cite{GN}.
In this latest paper the authors construct two particular approximate solutions of $Airy = 0$, denoted by $\phi^{app}_{f,\pm}$.
In particular, near $z = 0$, 
$$
\phi^{app}_{f,-} = {1 \over g'(z)} Ai(\gamma g(z))
$$
where $g(z)$ is Langer's transformation, which leads to 
$$
\phi^{app}_{f,-}(z) = Ai(\gamma (z - z_c)) + O(\nu^{1/4}) .
$$
The corresponding formulas for $\phi^{app}_{f,+}$ are similar, with $Ai$ replaced by $Ci$, a linear combination of $Ai$ and $Bi$.

We now construct an approximate Green's function for the $Airy$ operator. Let 
$$
G^{Ai}(z,y) = {1 \over \eps W^{Ai}(x)}  \left\{ 
\begin{array}{c}{\phi^{app}_{f,+}(z) \phi^{app}_{f,-}(x)} \quad \hbox{if} \quad z < x,
\\
 {\phi^{app}_{f,-}(z)  \phi^{app}_{f,+}(x)} \quad \hbox{if} \quad z > x,
\end{array}
\right.
$$
where $W^{Ai}$ is the Wronskian determinant of $\phi^{app}_\pm(x)$. 
Note that the Wronskian determinant is independent of $x$, since there is no first derivative term in $Airy$.
In addition, we have 
$$
W^{Ai}(x) \sim \gamma = O(\nu^{-1/4}).
$$ 
Therefore $G^{Ai}$ is rapidly decreasing in $z$ on both sides of $x$, within scales of order $\nu^{1/4}$.
By construction, 
$$
Airy \, G^{Ai}(x,z) = \delta_x + O(\nu^{3/4}) G^{Ai}(x,z) .
$$
Let us now turn to the study of the operator ${\cal A}$. Note that there exist four approximate independent solution to ${\cal A} = 0$.
Two are simply $\psi_{1,2} = \psi^{app}_{f,\pm}$ (which are approximate solutions of $Airy = 0$), and two others are defined by
$$
Airy(\psi_3) = e^{-\alpha z},  \qquad Airy(\psi_4) = {e^{\alpha z} - e^{- \alpha z} \over 2 \alpha} . 
$$
Note that $\psi_4$ does not decay to $0$. On the other side $\psi_3$ decays slowly
and will be pivotal in the relation dispersion for the adjoint operator. We explicitly have
$$
\psi_3 = G^{Ai} \star e^{- \alpha z} + O(\nu^{3/4}) .
$$
Note that $\psi_3$ is explicitly given by 
$$
\psi_3(x) = {1 \over \eps W^{Ai}} \Bigl( \phi_{f,-}^{app}(x)  \int_0^x \phi_{f,+}^{app}(y) e^{- \alpha y} dy
+ \phi_{f,+}^{app}(x) \int_x^{+ \infty} \phi_{f,-}^{app}(y) e^{- \alpha y} dy \Bigr) .
$$ 
Note also that for bounded $z$, and in particular in the critical layer $e^{-\alpha z}$ may replaced by $1$, up to $O(\alpha)$ 
smaller terms.

We define $\phi^{t,app}_{f,-} = \phi^{app}_{f,-}$, which is a solution of ${\cal A} = 0$ and is an approximate fast decaying mode.
Note that
\beq \label{partdisper1}
{\partial_z \phi^{app}_{f,-}(0) \over \phi^{app}_{f,-}(0)} = \gamma {Ai'(- \gamma z_c) \over Ai(- \gamma z_c)} + O(\nu^{1/4}) .
\eeq


\subsection{Construction of $\phi^{t,app}_{s,-}$}


The construction of $\phi^{t,app}_{s,-}$ is done by iteration.
We start from $\phi^t_{-,0} = 1$, which leads to a first guess 
$$
\psi_0 = e^{- \alpha z}.
$$
Note that the corresponding error is
\beq \label{e0}
e_0 = Orr^t(\psi_0) = - 2 \alpha U_s' e^{-\alpha z},
\eeq
which is of order $O(\alpha)$.
We thus introduce a first corrector to $\psi_0$ defined by
$$
Ray^t_\alpha(\psi_1) = - e_0.
$$ 
We then have
\beq \label{iteration2}
Orr^t(\psi_0 + \psi_1) = - \eps Diff(\psi_1),
\eeq
where $Diff(\psi_1)$ is defined in Section \ref{sec23}. We have
$$
\psi_1 = - (U_s(z) - c)^{-1} f_1,
$$
where
$$
f_1 = RaySolver_\alpha \Bigl( (U_s - c) e_0(x) \Bigr) = - 2 \alpha G_{R,\alpha} \star (U_s - c) U_s' e^{-\alpha z},
$$
where $G_{R,\alpha}$ is the Green function of Rayleigh's operator (see \cite{GN}).
In order to study the singularity of $\psi_1$ we write $f_1$ in the form of
$$
f_1 = f_1(z_c) + (U_s(z) - c) g_1(z),
$$
where $g_1$ is a smooth function. In particular $Diff(g_1)$ is also smooth.
Note that, as $e_0$ is of order $O(\alpha)$, $f_1(z_c)$ is also of order $O(\alpha)$.
Therefore $\psi_1$ has a singularity at $z_c$, which is exactly of the form 
$$
\psi_1 = \psi_1^1 + \psi_1^2 :=  - {f_1(z_c)  \over U_s (z) - c} - g_1(z).
$$ 
In particular, the leading term in  $-\eps Diff(\psi_1)$ is
$$
-\eps Diff(\psi_1) =  \eps f_1(z_c) Diff \Bigl( {1 \over U_s(z) - c}  \Bigr) + O(\eps). 
$$ 
To remove the singularity at leading order we introduce $\psi_2$ defined by
$$
\partial^2_z Airy \, \psi_2 =  \eps \partial_z^4 \psi_1^1 = - \eps f_1(z_c)  \partial_z^4 \Bigl ({ 1 \over U_s(z) - c} \Bigr).
$$
This leads to
\beq \label{theta2}
Airy \, \psi_2 =  -  \eps  f_1(z_c) \partial_z^2 \Bigl( { 1 \over U_s(z) - c } \Bigr).
\eeq
We recall that
$$
Airy \, \psi_2 = \Bigl( (U_s(z) -  c   ) - \eps \partial_z^2  \Bigr) \psi_2.
$$ 
Next we define $\psi_2^1$ by
$$
\psi_2 =    {f_1(z_c) \over U_s(z) -c }  + \psi_2^1.
$$
We have at leading order
$$
Airy \,   \psi_2^1 = - f_1(z_c),
$$
thus at leading order
$$
\psi_2^1 = - f_1(z_c) \psi_3(z),
$$
leading to
$$
\psi_2 = {f_1(z_c) \over U_s(z) -c } - f_1(z_c) \psi_3(z).
$$
Summing up, this gives an approximate mode of the form
\beq \label{approxmode}
\phi^{t,app}_{s,-}(z) = e^{-\alpha z} - f_1(z_c) \psi_3(z) - g_1(z) +\cdots.
\eeq
Note in particular that
\beq \label{disperstwo}
{\partial_z \phi^{t,app}_{s,-}(0) \over \phi^{t,app}_{s,-}(0) } = O(\nu^{-1/4}).
\eeq


\subsection{Dispersion relation}


The relation dispersion is
\beq \label{adjointdispersion}
{\partial_z \phi^{t,app}_{s,-}(0) \over \phi^{t,app}_{s,-}(0) }  = \gamma {Ai'(- \gamma z_c) \over Ai(- \gamma z_c)},
\eeq
which both are of size $O(\nu^{-1/4})$.
This dispersion relation should be the same as that of Orr Sommerfeld, namely (\ref{disper}), which is not
clear (and seems not proven, even at a formal level, in the physical literature).


\subsection{Use of the adjoint}


Let us recall in a general setting how to use the adjoint in order to study the perturbation of the spectrum of an operator. 
 Let $A(\eps)$ be a smooth families of linear operators, such that $A(\eps)$ has an eigenvector $u(\eps)$
with corresponding eigenvalue $\lambda(\eps)$. Assuming that $A(\eps)$, $u(\eps)$ and $\lambda(\eps)$ have expansions in $\eps$,
of the form 
$$
A = A_0 + \eps A_1 + \eps^2 A_2 + \cdots
$$ 
and similarly for $u$ and $\lambda$ and writing $A(\eps) = \lambda(\eps) u(\eps)$
we get
$$
A_0 e_0 = \lambda_0 e_0
$$ 
and at first order
\beq \label{first1}
(A_0 - \lambda_0) e_1 = \lambda_1 e_0 - A_1 e_0.
\eeq
Note that $(A_0 - \lambda_0)$ is precisely not invertible. Therefore we must choose $\lambda_1$ such that 
$$
\lambda_1 e_0 - A_1 e_0 \in Range(A_0 - \lambda_0) = ker(A_0^t - \lambda_0)^{\perp}.
$$ 
In particular
for any $v \in ker(A_0^t - \lambda_0)$,
$$
\lambda_1 (e_0 | v) = (A_1 e_0 | v) ,
$$
which gives $\lambda_1$. For this $\lambda_1$, (\ref{first1}) is precisely invertible, since the right hand side is in the range of $A_0 - \lambda_0$.
The construction may then be iterated to get a whole expansion of $\lambda(\eps)$ and $u(\eps)$.
The previous arguments are of course purely formal.

Conversely, if we know that $A(\eps)$ has such an expansion in powers of $\eps$, then, for any arbitrarily large $N$
we may construct approximate eigenvalues and normalized eigenvectors such that
$$
A(\eps) u(\eps) = \lambda(\eps) u(\eps) + O(\eps^N).
$$
Note that the adjoint is precisely used to construct an inverse of $(A_0 - \lambda_0)$ when $\lambda_0$ 
is in the spectrum of $A_0$.


\section{Other boundary conditions for Navier Stokes equations}


The previous linear analysis may be extended to nearby systems, namely to genuine Navier Stokes equations
with nearby boundary condition like succion or Navier boundary condition.
Let us detail the case of Navier boundary condition and let us consider  incompressible Navier Stokes equations in an half plane
\beq \label{NS1-navier-boundary} 
\partial_t u^\nu + (u^\nu  \cdot \nabla) u^\nu - \nu \Delta u^\nu + \nabla p^\nu = f^\nu,
\eeq
\beq \label{NS2-navier-boundary}
\nabla \cdot u^\nu = 0 ,
\eeq
together with
Navier boundary condition
\beq \label{navier}
u_1 = \beta \omega, \qquad u_2 = 0  \qquad \hbox{for} \qquad y = 0,
\eeq
where $u_1$ is the tangential velocity, $u_2$ the normal velocity, $\beta$ a constant and 
$\omega = -(\partial_z^2 \psi - \alpha^2 \psi)$. (Note that $\omega=\nabla\times u=\partial_xu_2-\partial_zu_1$.)

We are interested in the stability of a shear layer profile
$U(z) = (U_s(z),0)$. The objective of this section is to prove that long wave instabilities persist provided $\beta$ is small enough.
More precisely we will prove the following claim.

\medskip

{\it Claim: There exists $\beta_0$, small enough, such that 
if $\beta \le \beta_0 \nu^{1/4}$ then there exists a linearly unstable mode to Navier Stokes equations with Navier boundary condition.}

\medskip

As for the Dirichlet boundary condition we look for velocities of the form
$$
u= \nabla^\perp \Bigl( e^{i \alpha (x - c t) } \psi(z) \Bigr)=e^{i \alpha (x - c t) } (\partial_z \psi, -i\alpha \psi),
$$
which also now gives to the same Orr Sommerfeld operator, except that the Navier boundary condition (\ref{navier})
leads to 
$$
\psi(0) =0, \qquad  \partial_z \psi(0) =- \beta \partial_z^2\psi(0).
$$
As previously we must find a linear combination of $\psi_{s,-}$ and $\psi_{f,-}$ which satisfies (\ref{navier}).
First note that, $\omega = \partial_z^2 \psi - \alpha^2 \psi$. Neglecting the second term, this gives
$$
a \psi_{f,-}(0) + b \psi_{s,-} (0) = 0,
$$
$$
a \partial_z \psi_{f,-}(0) + b \partial_z \psi_{s,-}(0) = -
\beta \Bigl( a \partial_z^2 \psi_{f,-}(0) + b \partial_z^2 \psi_{s,-}(0)  \Bigr) .
$$
Note that at leading order
$$
\partial_z^2 \psi_{f,-}(0) = \gamma^2 Ai(-\gamma z_c) .
$$
This leads to, neglecting $\beta \partial_z^2 \psi_{s,-}(0)$ in front of $\partial_z \psi_{s,-}(0)$
\beq \label{disperbeta}
{   \psi_{f,-} (0) \over \partial_z  \psi_{f,-} (0) +\beta \partial_z^2 \psi_{f,-}(0)}
= {  \psi_{s,-} (0) \over \partial_z  \psi_{s,-} (0)}
\eeq
and to the dispersion relation
\beq \label{disperbeta2}
\alpha {U_+^2 \over U_s'(0)^2}  - {c \over U_s'(0)} = \gamma^{-1} {Ai(2, - \gamma z_c) \over 
Ai(1,-\gamma  z_c) +\beta \gamma Ai(- \gamma z_c)} 
+O(\alpha^2).
\eeq
The study of the dispersion relation continues as in the previous section, except that the Tietjens' function is replaced
by 
$$
Ti(z) = {Ai(2,z e^{- 5 i \pi / 6}) \over z e^{- 5 i \pi / 6} \Bigl( Ai(1,z e^{- 5 i \pi / 6})  
+i^{1/3} \alpha_0^{1/3}  \beta \nu^{-1/4} U_s'(z_c)^{1/3} Ai(z e^{-5 i \pi / 6})  \Bigr)} .
$$
We therefore look at $\beta$ of order $\nu^{1/4}$ and define
$$
\beta = \beta_0 \nu^{1/4} .
$$
A continuity argument immediately gives that provided $\beta_0$ is small enough, there exists an unstable mode
to Navier Stokes equations with Navier boundary condition, which ends the proof of the claim.

\section{Rotation or magnetic field}


In this section we focus on Navier Stokes equations in a rotating frame, 
which read
 \beq \label{NS1C} 
\partial_t u^\nu + (u^\nu  \cdot \nabla) u^\nu - \nu \Delta u^\nu + \eta e \times u^\nu + \nabla p^\nu = f^\nu,
\eeq
\beq \label{NS2C}
\nabla \cdot u^\nu = 0 ,
\eeq
where $e$ is a fixed vector and $\eta$ a scalar. Note that $\eta e \times u$ is the usual Coriolis force
(or magnetic force in the case of a charged fluid).
We will focus on the case where $e$ is perpendicular to the plane $z = 0$.
The case of an arbitrarily $e$ can be treated using similar ideas.

We will prove the following claim.

\medskip

{\it Claim: there exists $N$ large enough ($N \ge 3$) such that if $\eta \le \nu^N$, then there exists an approximate
linear unstable mode to \eqref{NS1C} and \eqref{NS2C}.}

\medskip

If $\eta \to + \infty$ as $\nu \to 0$, Ekman's layers appear on the boundary, with a typical size
$\sqrt{\nu / \eta}$. These layers are known to be linearly
stable provided the associated Reynolds number  $Re = \| u^\nu \|_{L^\infty} / \sqrt{\eta \nu}$ remains small enough. 
If $\eta^{-1}$ is of order $\nu$, they are known to be linearly and nonlinear unstable for this Reynolds number
 large enough.

In this paper we are interested in the case where $\eta$ is small. In this case the size of the boundary layer is
$O(\sqrt{\nu})$, as for genuine Prandtl's layer. In this paper we focus on the evolution of the stability of such
layers as $\eta$ increases, but remains small.
The case where $\eta$ is large, but negligible with respect to $\nu^{-1}$ remains open.
 
We start with a shear layer, as for Navier Stokes equations. We first note that the Coriolis force $e \times u$ may be 
absorbed into the pressure term.

For sufficiently small $\eta$ we construct a genuine three dimensional instability starting from a two dimensional one.
Note that one of the main interest of this construction is that the instability is fully three dimensional, and not two dimensional
as for genuine Navier Stokes equations.
Namely let us choose $e_1$, direction of the $x$ variable, parallel to the plane $ z = 0$. 
Using Orr Sommerfeld equation, we may construct
an instability $u_0$ which is independent on $y$ and only depends on $x$ and $z$.
This instability $u_0$ solves (\ref{NS1C}) up to the Coriolis term $\eta e \times u_0$, a term in the direction $y$.
We therefore alter Orr Sommerfeld to take into account this velocity in the $y$ direction. More precisely we look for
$(u,v,w)$ of the form
$$
(u,v,w) = e^{i \alpha (x - c t)} ( \psi'(z), v(z), -i \alpha \psi(z) ).
$$
This leads to the following modified Orr Sommerfeld equations
\beq \label{OrrCor1}
 (U_s - c)  (\partial_z^2 - \alpha^2) \psi - U_s''  \psi  - { \nu \over i \alpha}   (\partial_z^2 - \alpha^2)^2 \psi 
 - {\eta \over i \alpha} \partial_z v = 0,
 \eeq
\beq \label{OrrCor2}
(U_s - c) v - \eps (\partial_z^2 - \alpha^2) v - {\eta \over i \alpha} \psi' = 0.
\eeq
Note that \eqref{OrrCor1}-\eqref{OrrCor2} is a coupling between a genuine Orr Sommerfeld equation and
an Airy type equation, through terms of order $\eta / i \alpha$.
Let
$$
Airy(v) = (U_s - c + \eps \alpha^2) v - \eps \partial_z^2 v.
$$
Let us study the evolution of an eigenvector $\psi(\eta)$, $v(\eta)$ with corresponding eigenvalue $c(\eta)$
as $\eta$ increases from $0$. We will only detail formal computations and look for an asymptotic
expansion of these quantities in $\eta$. This leads to the construction of approximate unstable eigenmodes,
a first step in the construction of nonlinear instabilities. We have, underlying the $c$ dependency of the operators on $c$,
$$
Orr(c,\psi) =  {\eta \over i \alpha} \partial_z v,
$$
$$
Airy(c,v) =  {\eta \over i \alpha} \psi'.
$$
We look for expansions of the form 
$$
c = c_0 + \eta^2 c_1 + \eta^4 c_2 + \cdots,
$$ 
$$
\psi = \psi_0 + \eta^2 \psi_1 + \eta^4 \psi_2 + \cdots
$$
and
$$
v = \eta v_0 + \eta^3 v_1 + \cdots.
$$
This leads to $Orr(c_0,\psi_0) = 0$, meaning that we construct a branch of solution starting from a long wave instability
$(c_0,\psi_0)$ of the non rotating Navier Stokes equations. Moreover, 
$$
Airy(c_0,v_0) = {1 \over i \alpha} \psi_0'.
$$
The next order gives 
\beq \label{nextorder}
Orr(c_0,\psi_1) + Orr(c_1,\psi_0) =  {1 \over i \alpha} \partial_z v_0 + c_0 (\partial_z^2 - \alpha^2) \psi_0.
\eeq
Let $\psi_1^t$ be the eigenvalue of $Orr^t$ with the same eigenvalue $c_0$. Then (\ref{nextorder}) may be solved
provided
$$
\Bigl( Orr(c_1,\psi_0), \psi_1^t \Bigr) =   \Bigl(-i \alpha^{-1} \partial_z v_0 + c_0 (\partial_z^2 - \alpha^2) \psi_0 , \psi_1^t \Bigr),
$$
or expanding the Orr Sommerfeld operator, and denoting by $v_1$, $\omega_1$, $v_1^t$ and $\omega_1^t$ the corresponding
velocity fields and vorticities we get 
\beq \label{c-1-solution}
\begin{split}
c_1 \int v_1 \bar v_1^t &=  -i \alpha^{-1} ( \partial_zv_0, \psi_1^t)
- \int U_s \omega_1 \bar \psi_1^t + \int U_s'' \psi_1 \bar \psi_1^t\\
& \quad+ {\nu \over i \alpha} \int \omega_1 \bar \omega_1^t+\int c_0 \omega_1 \bar \psi_1^t,
\end{split}
\eeq
which explicitly gives $c_1$. The formal construction of an approximate eigenmode can then be continued, up
to any arbitrarily high order. Keeping track of the various dependencies in $\nu$ we get that, provided $\eta \le \nu^N$
with $N$ large enough ($N \ge 3$), their exists an approximate unstable eigenmode to \eqref{NS1C} and \eqref{NS2C}.
Using the arguments developed in \cite{GN1} it is then possible to get an exact description of the corresponding
Green function, and an exact unstable eigenmode. We will not detail these points here, which are routine, but lengthly.


\section{Effect of density}


Similar technics may be used to study compressible Navier Stokes equations or inhomogeneous Navier Stokes equations,
provided the effects of compressibility or inhomogeneity are very small. Let us detail the case of compressible fluids.

Let us focus on the classical barotropic compressible Navier Stokes equations
\beq \label{NSC1}
\partial_t \rho + \nabla  \cdot (\rho u) = 0,
\eeq
\beq \label{NSC2}
\partial_t  u + (u  \cdot\nabla) u - \mu \rho^{-1}  \Delta u - \xi \rho^{-1} \nabla (\nabla  \cdot u) + {\nabla p(\rho) \over \eta^2 \rho} = 0
\eeq
with boundary conditions $u = 0$ on the boundary $y = 0$.
In these equations, $\eta$ is the Mach number, which is assumed to be small, and $\rho$ is a small
perturbation of a constant, fixed to be $1$. As $\eta$ goes to $0$, formally $\rho$ goes to $0$, therefore
$\nabla \cdot u$ goes to $0$ and we recover incompressible Navier Stokes equations.
Let $h(\rho)$ such that
$$
\nabla h(\rho) = {\nabla p(\rho) \over \rho} .
$$
The objective of this paper is to study the stability of the shear layer profile $\rho_s = 1$, 
$U_s = (U_s(y),0)$. More precisely we will prove formally that such a shear profile exhibits long wave instabilities
provided $\eta$ is small enough.

The linearized system reads
\beq \label{NSC1lin}
\partial_t \rho + \nabla \cdot (u + \rho U_s) = 0,
\eeq
\beq \label{NSC2lin}
\partial_t  u + ( U_s  \cdot \nabla) u  + ( u \cdot \nabla) U_s 
- \mu \Delta u  + \mu \rho \Delta U_s - \xi \nabla (\nabla  \cdot u) + h'(1) {\nabla \rho \over \eta^2} = 0.
\eeq
We will prove the following claim.

\medskip

{\it Claim: there exists $N$ large enough such that if $\eta \le \nu^N$, then there exists a linear unstable mode to
\eqref{NSC1lin} and \eqref{NSC2lin}.}

\medskip

The strategy is to construct an approximate growing mode of \eqref{NSC1lin} and \eqref{NSC2lin}  starting from 
a long wave growing mode of the genuine incompressible Navier Stokes equations.
As usual in such incompressible limits, $\rho$ is of order $O(\eta^2)$. Moreover, the velocity
field will only be divergence free at leading order. We will therefore look for velocities and density of
the form
$$
(u,w) = (u_0,w_0) + \eta^2 (u_2,w_2) + \eta^2 \nabla \theta_2 +  \cdots,
$$
$$
\rho = \eta^2 \rho_2 + \cdots
$$
and
$$
c = c_0 + \eta^2 c_2 + \cdots .
$$
We start with a linear instability of incompressible Navier Stokes equations, of the form
$$
(u_0,w_0) = e^{i \alpha (x - c t)} ( \psi'(z), -i \alpha \psi(z) ) .
$$
This instability has a corresponding pressure $p_0$. As $\rho$ will be of order $O(\eta^2)$, $\mu \rho \Delta U_s$
is negligible at this stage, and  moreover $\nabla \cdot u = 0$, hence we define $\rho_2$ by
$$
h'(1) \nabla \rho_2 = - \nabla p_0 .
$$
Now using (\ref{NSC1lin}) we get
$$
\Delta \theta_2 = - \partial_t \rho_2 - \nabla  \cdot (\rho_2 U_s),
$$
which defines $\theta_2$. We then project (\ref{NSC2lin}) on divergence free vector fields and gradient vector fields.
The project on gradients allows to define $\rho_4$. The projection on divergence free vector fields gives $(u_2,w_2)$,
up to a solvability condition which gives $c_2$. This construction can be iterated at any order.
This defines an approximate unstable mode. Using technics of \cite{GN1} we can study the corresponding Green functions
and construct an exact unstable mode. This proves the claim, by keeping track of the $\nu$ dependency of the various operators
(we do not detail the precise computation of $N$).

\bigskip
\noindent {\bf Acknowledgment.}
D. Bian is supported by NSFC under the contract 11871005. 
\small



\end{document}